\documentclass{gtart_h}

\def\ifplaintex{\expandafter\ifx\csname documentclass\endcsname\relax}

\def\gtp{{\mathsurround=0pt\it $\cal G\mskip-2mu$eometry \&\ 
$\cal T\!\!$opology $\cal P\!$ublications}}  

\def\recd{{\small Received:\qua\receiveddate\ifx\reviseddate\relax
\else\qquad Revised:\qua\reviseddate\fi\par}} 

\def\lognumber#1{\def\thelognumber{#1}}
\def\volumenumber#1{\def\thevolumenumber{#1}}
\def\volumeyear#1{\def\thevolumeyear{#1}}
\def\papernumber#1{\def\thepapernumber{#1}}
\def\pagenumbers#1#2{\def\startpage{#1}\def\finishpage{#2}}
\def\published#1{\def\publishdate{#1}}

\def\received#1{\def\receiveddate{#1}}
\def\revised#1{\def\reviseddate{#1}}
\def\accepted#1{\def\accepteddate{#1}}

\def\asciiaddress#1{\def\theasciiaddress{#1}}

\long\def\asciiabstract#1{\long\def\theasciiabstract{#1}}

\let\\\par\let\thelognumber\relax\let\thevolumenumber\relax
\let\thepapernumber\relax\let\thevolumeyear\relax\let\startpage\relax
\let\finishpage\relax\let\publishdate\relax\let\receiveddate\relax
\let\reviseddate\relax\let\accepteddate\relax\let\theasciititle\relax
\let\theasciiauthors\relax\let\theasciiaddress\relax
\let\theasciiabstract\relax

\let\theasciiemail\relax

\ifplaintex
\font\logobig=cmssbx10 scaled 3836
\font\logomed=cmssbx10 scaled 2557
\else
\font\logobig=cmssbx10 scaled 4200
\font\logomed=cmssbx10 scaled 2800
\fi

\long\def\makeagttitle{   
\count0=\startpage
\agt\hfill      
\hbox to 45truept{\vbox to 0pt{\vglue -13truept{\logomed A\kern -.37em{\logobig 
T}\kern -.38em G}\vss}\hss}
\break
{\small Volume \thevolumenumber\ (\thevolumeyear)
\startpage--\finishpage\nl
Published: \publishdate}

\vglue .25truein

{\parskip=0pt\leftskip 0pt plus
1fil\def\\{\par\smallskip}{\Large\bf\thetitle}\par\medskip} \vglue
0.05truein

{\parskip=0pt\leftskip 0pt plus 1fil\def\\{\par}{\sc\theauthors}
\par\medskip}%
 
\vglue 0.03truein 

{\small\leftskip 25truept\rightskip 25truept{\bf Abstract}\stdspace\theabstract

{\bf AMS Classification}\stdspace\theprimaryclass
\ifx\thesecondaryclass\relax\else; \thesecondaryclass\fi\par
{\bf Keywords}\stdspace \thekeywords\par}\vglue 7truept

}   

\ifplaintex
\font\phead=cmsl9 scaled 950
\font\pnum=cmbx10 scaled 913
\font\pfoot=cmsl9 scaled 950
\headline{\vbox to 0pt{\vskip -4.5mm\line{\small\phead\ifnum
\count0=\startpage ISSN 1472-2739 (on-line) 1472-2747 (printed)
\hfill {\pnum\folio}\else\ifodd\count0\def\\{ }%
\ifx\theshorttitle\relax\thetitle\else\theshorttitle\fi\hfill{\pnum\folio}
\else\def\\{ and }{\pnum\folio}\hfill\ifx\theshortauthors\relax\theauthors
\else\theshortauthors\fi\fi\fi}\vss}}
\footline{\vbox to 0pt{\vglue 0mm\line{\small\pfoot\ifnum\count0=\startpage
\copyright\ \gtp\hfill\else
\agt, Volume \thevolumenumber\ (\thevolumeyear)\hfill\fi}\vss}}
\else
\font=cmsl9 scaled 1050
\font\lnum=cmbx10 
\font=cmsl9 scaled 1050
\makeatletter
\def\@oddhead{{\small\lhead\ifnum\count0=\startpage ISSN 1472-2739 
(on-line) 1472-2747 (printed)\hfill {\lnum\number\count0}\else\ifodd\count0
\def\\{ }\ifx\theshorttitle\relax \thetitle \else\theshorttitle\fi\hfill
{\lnum\number\count0}\else\def\\{ and }{\lnum\number\count0}
\hfill\ifx\theshortauthors\relax 
\theauthors\else\theshortauthors\fi\fi\fi}}\def\@evenhead{\@oddhead}
\def\@oddfoot{\small\lfoot\ifnum\count0=\startpage\copyright\ \gtp\hfill\else
\agt, Volume \thevolumenumber\ (\thevolumeyear)\hfill\fi}
\def\@evenfoot{\@oddfoot}
\makeatother
\fi
\let\maketitlepage\makeagttitle
\let\makeshorttitle\maketitlepage
\let\maketitle\maketitlepage

\newwrite\gtoutfile
\long\gdef\makeheadfile{  
{\def\\{, }\def\s{ }
\immediate\openout\gtoutfile head.xxx
\immediate\write\gtoutfile{Proxy-for: \ifx\theasciiauthors\relax
\theauthors\else\theasciiauthors\fi\s<\ifx\theasciiemail\relax\theemail\else\theasciiemail\fi>}
\immediate\write\gtoutfile{\noexpand\\}
\immediate\write\gtoutfile{Authors: \ifx\theasciiauthors\relax
\theauthors\else\theasciiauthors\fi}
{\def\\{ }\immediate\write\gtoutfile{Title: \ifx\theasciititle\relax
\thetitle\else\theasciititle\fi}}
\immediate\write\gtoutfile{Subj-class: GT or SG, GR etc}
\immediate\write\gtoutfile{MSC-class: \theprimaryclass\ifx\thesecondaryclass\relax\else, \thesecondaryclass\fi}
\immediate\write\gtoutfile{Journal-ref: Algebr. Geom. Topol. \thevolumenumber\s
(\thevolumeyear) \startpage-\finishpage}
\immediate\write\gtoutfile{Comments: Published by Algebraic and
Geometric Topology at}
\immediate\write\gtoutfile{\s\s\s  http://www.maths.warwick.ac.uk/agt/AGTVol\thevolumenumber/agt-\thevolumenumber-\thepapernumber.abs.html}
\immediate\write\gtoutfile{\noexpand\\}
\immediate\write\gtoutfile{}
\ifx\theasciiabstract\relax
\immediate\write\gtoutfile{\theabstract}\else
\immediate\write\gtoutfile{\theasciiabstract}\fi
\immediate\write\gtoutfile{}
\immediate\write\gtoutfile{\noexpand\\}
\immediate\write\gtoutfile{}
\immediate\closeout\gtoutfile}}  

\def\maketitlepage{\makeagttitle\makeheadfile}
\let\makeshorttitle\maketitlepage
\let\maketitle\maketitlepage

\lognumber{23}
\volumenumber{5}
\volumeyear{2005}
\papernumber{23}
\pagenumbers{537}{562}
\received{16 September 2004} 
\revised{18 May 2005}
\accepted{3 June 2005}
\published{19 June 2005}

\usepackage{amsmath,amssymb,graphicx}

\def\textsection#1{subsection #1}

\theoremstyle{plain}
  \newtheorem{theorem}{Theorem}
  \newtheorem{lemma}[theorem]{Lemma}
  \newtheorem{proposition}[theorem]{Proposition}
  \newtheorem{corollary}[theorem]{Corollary}
\theoremstyle{definition}
  \newtheorem{definition}[theorem]{Definition}
  \newtheorem{question}[theorem]{Question}
  \newtheorem{remark}[theorem]{Remark}
  \newtheorem{example}[theorem]{Example}

\newcommand{\Z}{\mathbb{Z}}                       
\newcommand{\Q}{\mathbb{Q}}                       
\newcommand{\K}{\mathbb{K}}                       
\newcommand{\A}{\mathbb{A}}                       
\newcommand{\m}{\mathfrak{m}}                     
\newcommand{\fps}[1]{\mathopen{[\![}#1\mathclose{]\!]}}   

\newcommand{\Ker}{\operatorname{ker}}             
\renewcommand{\Im}{\operatorname{im}}             
\newcommand{\Hom}{\operatorname{Hom}}             
\newcommand{\End}{\operatorname{End}}             
\newcommand{\Aut}{\operatorname{Aut}}             
\newcommand{\Inn}{\operatorname{Inn}}             
\newcommand{\Br}[1]{\operatorname{B}_{#1}}        

\newcommand{\tensor}[1][]{\mathbin{\otimes_{#1}}} 
\newcommand{\compose}{\mathbin{\circ}}            
\newcommand{\tr}{\operatorname{tr}}               
\newcommand{\id}{\operatorname{id}}               
\newcommand{\0}{\cdot}                            
\newcommand{\I}{\makebox{\rm I}}                  
\newcommand{\II}%
{\makebox{\makebox[1.2pt][l]{\I}{\I}}}            

\newcommand{\HYB}{H_{\scriptscriptstyle\rm YB}}   
\newcommand{\indices}[2]%
{{\left[\begin{matrix}#1\\#2\end{matrix}\right]}} 
\newcommand{\smallindices}[2]%
{{\left[\begin{smallmatrix}#1\\#2\end{smallmatrix}\right]}} 

\begin{document} 

\title{Yang-Baxter deformations of quandles and racks}

\author{Michael Eisermann}

\address{Institut Fourier, Universit\'e Grenoble I, 38402 St Martin d'H\`eres, 
France}
\asciiaddress{Institut Fourier, Universite Grenoble I, 38402 St Martin d'Heres, 
France}

\email{Michael.Eisermann@ujf-grenoble.fr}

\urladdr{www-fourier.ujf-grenoble.fr/~eiserm}

\begin{abstract}
Given a rack $Q$ and a ring $\A$, one can construct 
a Yang-Baxter operator $c_Q\colon V\tensor V\to V\tensor V$ 
on the free $\A$-module $V = \A{Q}$ by setting 
$c_Q(x \tensor y) = y \tensor x^y$ for all $x,y\in Q$. 
In answer to a question initiated by D.N.\,Yetter and P.J.\,Freyd, this article 
classifies formal deformations of $c_Q$ in the space of Yang-Baxter operators.
For the trivial rack, where $x^y = x$ for all $x,y$, one has, of course, 
the classical setting of r-matrices and quantum groups.
In the general case we introduce and calculate the cohomology theory 
that classifies infinitesimal deformations of $c_Q$. 
In many cases this allows us to conclude that $c_Q$ is rigid.
In the remaining cases, where infinitesimal deformations are possible, 
we show that higher-order obstructions 
are the same as in the quantum case.
\end{abstract}

\asciiabstract{Given a rack Q and a ring A, one can construct a
Yang-Baxter operator c_Q: V tensor V --> V tensor V  on the free
A-module V = AQ by setting c_Q(x tensor y) = y tensor x^y for all x,y
in Q.  In answer to a question initiated by D.N. Yetter and
P.J. Freyd, this article classifies formal deformations of c_Q in the
space of Yang-Baxter operators.  For the trivial rack, where x^y =
x for all x,y, one has, of course, the classical setting of r-matrices
and quantum groups.  In the general case we introduce and calculate
the cohomology theory that classifies infinitesimal deformations of
c_Q.  In many cases this allows us to conclude that c_Q is rigid.  In
the remaining cases, where infinitesimal deformations are possible, we
show that higher-order obstructions are the same as in the quantum
case.}

\keywords{Yang-Baxter operator, r-matrix, braid group representation,
deformation theory, infinitesimal deformation, Yang-Baxter cohomology}

\primaryclass{17B37}
\secondaryclass{18D10,20F36,20G42,57M25}

\maketitle

\section*{Introduction}

Following M.\,Gerstenhaber \cite{Gerstenhaber:1964}, 
an algebraic deformation theory should 
\begin{itemize}
\item
define the class of objects within which deformation takes place,
\item
identify infinitesimal deformations as elements of a suitable cohomology,
\item
identify the obstructions to integration of an infinitesimal deformation,
\item
give criteria for rigidity, and possibly determine the rigid objects.
\end{itemize}

In answer to a question initiated by P.J.\,Freyd 
and D.N.\,Yetter \cite{FreydYetter:1989}, we carry out 
this programme for racks (linearized over some ring $\A$) and 
their formal deformations in the space of $\A$-linear Yang-Baxter operators.


A rack is a set $Q$ with a binary operation, denoted $(x,y) \mapsto x^y$, 
such that $c_Q\colon x \tensor y \mapsto y \tensor x^y$ defines 
a Yang-Baxter operator on the free $\A$-module $\A{Q}$ 
(see Section \ref{sec:BasicsAndResults} for definitions). 
For a trivial rack, where $x^y = x$ for all $x,y\in Q$, 
we see that $c_Q$ is simply the transposition operator.
In this case the theory of quantum groups 
\cite{Drinfeld:1987,Turaev:1988,Kassel:1995,KasselRossoTuraev:1997}
produces a plethora of interesting deformations, 
which have received much attention over the last 20 years.
It thus seems natural to study deformations of $c_Q$ 
in the general case, where $Q$ is a non-trivial rack.

\subsection*{Outline of results}

We first introduce and calculate the cohomology theory 
that classifies infinitesimal deformations of racks 
in the space of Yang-Baxter operators.
In many cases this suffices to deduce rigidity.
In the remaining cases, where infinitesimal deformations are possible, 
we show that higher-order obstructions do not depend on $Q$:
they are the same as in the classical case of quantum invariants.
(See \textsection\ref{sub:Deformations} for a precise statement.)

Formal Yang-Baxter deformations of racks 
thus have an unexpectedly simple description: 
up to equivalence they are just r-matrices with a special symmetry 
imposed by the inner automorphism group of the rack.
Although this is intuitively plausible, it requires 
a careful analysis to arrive at an accurate formulation.
The precise notion of \emph{entropic} r-matrices will 
be defined in \textsection\ref{sub:EntropicMaps}. 

With regards to topological applications, this result may 
come as a disappointment in the quest for new knot invariants. 
To our consolation, we obtain a complete and concise solution 
to the deformation problem for racks, which is quite 
satisfactory from an algebraic point of view.

Throughout our calculations we consider the generic case
where the order $|\Inn(Q)|$ of the inner automorphism group 
of the rack $Q$ is invertible in the ground ring $\A$.
We should point out, however, that certain knot invariants arise 
only in the modular case, where $|\Inn(Q)|$ vanishes in $\A$; 
see the closing remarks in Section \ref{sec:OpenQuestions}.

\subsection*{How this paper is organized}

In order to state the results precisely, 
and to make this article self-contained,
Section \ref{sec:BasicsAndResults} first recalls the notions of 
Yang-Baxter operators (\textsection\ref{sub:YangBaxterOperators}) 
and racks (\textsection\ref{sub:Racks}).
We can then introduce entropic maps (\textsection\ref{sub:EntropicMaps})
and state our results (\textsection\ref{sub:Deformations}).
We also discuss some examples (\textsection\ref{sub:Examples})
and put our results into perspective by briefly reviewing 
related work (\textsection\ref{sub:RelatedWork}).

The proofs are given in the next four sections:
Section \ref{sec:YangBaxterCohomology} introduces Yang-Baxter cohomology
and explains how it classifies infinitesimal deformations.
Section \ref{sec:RackCohomology} calculates 
this cohomology for racks.
Section \ref{sec:CompleteDeformations} generalizes 
our classification from infinitesimal to complete deformations.
Section \ref{sec:HigherObstructions} examines 
higher-order obstructions and shows that they are 
the same as in the classical case of quantum invariants.
Section \ref{sec:OpenQuestions}, finally, discusses some open questions.


\section{Review of basic notions and statement of results}
\label{sec:BasicsAndResults}

\subsection{Yang-Baxter operators} 
\label{sub:YangBaxterOperators}

Let $\A$ be a commutative ring with unit.
In the sequel all modules will be $\A$-modules, 
and all tensor products will be formed over $\A$.
For every $\A$-module $V$ we denote by $V^{\tensor n}$ 
the $n$-fold tensor product of $V$ with itself.
The identity map of $V$ is denoted by $\I\colon V\to V$,
and $\II = \I\tensor\I$ stands for the identity map of $V \tensor V$.

\begin{definition}
A \emph{Yang-Baxter operator} on $V$ is an automorphism 
$c\colon V \tensor V \to V \tensor V$ that satisfies 
the \emph{Yang-Baxter equation}, also called \emph{braid relation},
$$
(c \tensor\I)(\I\tensor c)(c \tensor\I) = 
(\I\tensor c)(c \tensor\I)(\I\tensor c) 
\qquad\text{in}\quad \Aut_\A(V^{\tensor3}).
$$
\end{definition}

\begin{figure}[hbtp]
  \centering
  \includegraphics{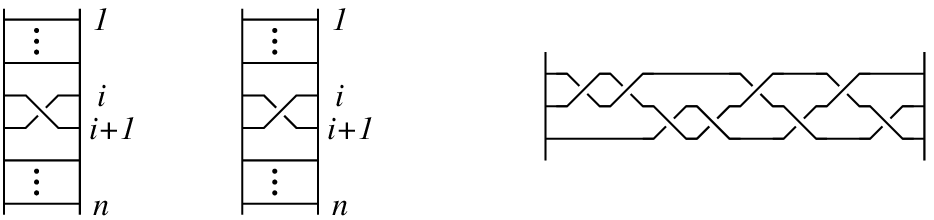}
  \caption{Elementary braids $\sigma_i^{+1}$, $\sigma_i^{-1}$; a more complex example $\beta$}
  \vspace{-1mm}
  \label{fig:Braids}
\end{figure}

This equation first appeared in theoretical 
physics, in a paper by C.N.\,Yang \cite{Yang:1967} 
on the many-body problem in one dimension,
in work of R.J.\,Baxter \cite{Baxter:1972,Baxter:1982}
on exactly solvable models in statistical mechanics,
and later in quantum field theory \cite{Faddeev:1984}
in connection with the quantum inverse scattering method.
It also has a very natural interpretation in terms of Artin's braid groups 
\cite{Artin:1947,Birman:1974} and their tensor representations:

\begin{remark}
Recall that the braid group on $n$ strands can be presented as
$$
\Br{n}  = \left\langle\; \sigma_1,\dots,\sigma_{n{-}1} \;\Big|\;
\begin{matrix} 
\hfill{} \sigma_i \sigma_j = \sigma_j \sigma_i  \hfill{} \quad \text{for } |i-j|\ge2 \\
\sigma_i \sigma_j \sigma_i  = \sigma_j  \sigma_i \sigma_j \quad \text{for } |i-j|=1
\end{matrix} \;\right\rangle,
$$
where the braid $\sigma_i$ performs a positive half-twist of the strands $i$ and $i+1$.
In graphical notation, braids can conveniently be represented as in Figure \ref{fig:Braids}.

Given a Yang-Baxter operator $c$, 
we can define automorphisms $c_i\colon V^{\tensor n} \to V^{\tensor n}$ for $i=1,\dots,n-1$
by setting $c_i = \I^{\tensor(i{-}1)} \tensor c \tensor \I^{\tensor(n{-}i{-}1)}$.
The Artin presentation ensures that there exists a unique braid group 
representation $\rho_c^n\colon \Br{n} \to \Aut_\A(V^{\tensor n})$
defined by $\rho_c^n(\sigma_i) = c_i$. 

Here we adopt the following convention:
braid groups will act on the left, so that composition 
of braids corresponds to composition of maps.
The braid in Figure \ref{fig:Braids}, for example, reads 
$\beta = \sigma_1^{-2}\sigma_2^2\sigma_1^{-1}\sigma_2^1\sigma_1^{-1}\sigma_2^1$;
it is represented by the operator $\rho_c^3(\beta) = 
c_1^{-2} c_2^2 c_1^{-1} c_2^1 c_1^{-1} c_2^1$ acting on $V^{\tensor 3}$.
\end{remark}

Notice that Artin, after having introduced his braid groups, 
could have written down the Yang-Baxter equation in the 1920s,
but without any non-trivial examples the theory would have remained void.
It is a remarkable fact that the Yang-Baxter equation 
admits any interesting solutions at all.
Many of them have only been discovered since the 1980s,
and our first example recalls the most prominent one:

\begin{example} \label{exm:Jones}
For every $\A$-module $V$ the transposition $\tau\colon V \tensor V \to V \tensor V$
given by $\tau(a \tensor b) = b \tensor a$ is a Yang-Baxter operator. This in itself
is not very surprising, but deformations of $\tau$ can be very interesting:
Suppose that $V$ is free of rank $2$ and choose a basis $(v,w)$. If we equip 
$V^{\tensor 2}$ with the basis $(v\tensor v, v\tensor w, w\tensor v, w\tensor w)$
then $\tau$ is represented by the matrix $c_1$ as follows:
$$
c_1 = 
\begin{pmatrix}
1 & 0 & 0 & 0 \\
0 & 0 & 1 & 0 \\
0 & 1 & 0 & 0 \\
0 & 0 & 0 & 1
\end{pmatrix}, \qquad
c_q =
\begin{pmatrix}
      q &       0 &       0 &     0 \\
      0 &       0 &     q^2 &     0 \\
      0 &     q^2 & q{-}q^3 &     0 \\
      0 &       0 &       0 &     q
\end{pmatrix}.
$$
For every choice of $q\in\A^\times$, the matrix $c_q$ is 
a Yang-Baxter operator, and for $q=1$ we obtain the initial solution $c_1 = \tau$. 
The family $(c_q)$, together with a suitable trace, yields 
the celebrated Jones polynomial \cite{Jones:1985,Jones:1987,Jones:1989}, 
a formerly unexpected invariant of knots and links.
More generally, deformations of $\tau$ lead 
to the so-called \emph{quantum invariants} of knots and links.
\end{example}

Given the matrix $c_q$ of the preceding example, it is straightforward
to \emph{check} that it satisfies the Yang-Baxter equation.
How to \emph{find} such solutions, however, is a much harder question.
Attempts to construct solutions in a systematic way have led 
to the theory of quantum groups (cf.\ \cite{Drinfeld:1987}).
For details we refer to the concise introduction 
\cite{KasselRossoTuraev:1997} or the textbook 
\cite{Kassel:1995}.

\begin{remark}
A slight reformulation sometimes proves useful.
Every Yang-Baxter operator $c$ can be written as $c = \tau f$
where $f \colon V \tensor V \to V \tensor V$ is an automorphism
satisfying $f_{12} f_{13} f_{23} = f_{23} f_{13} f_{12}$, with $f_{ij}$ 
acting on the $i$th and $j$th factor of $V \tensor V \tensor V$.
Such an operator $f$ is called an \emph{r-matrix}.
Depending on the context it may be more convenient 
to consider the r-matrix $f$ or the Yang-Baxter operator $\tau f$.
\end{remark}

The set of Yang-Baxter operators is closed under conjugation by $\Aut(V)$,
and conjugate operators yield conjugate braid group representations.
A general goal of Yang-Baxter theory, as yet out of reach, would be 
to classify solutions of the Yang-Baxter equation modulo conjugation by $\Aut(V)$.
In favourable cases this can be done at least \emph{locally}, that is, 
one can classify \emph{deformations} of a given Yang-Baxter operator.
Our main result, as stated in \textsection\ref{sub:Deformations} below,
covers a large class of such examples.

\begin{definition}
We fix an ideal $\m$ in the ring $\A$.
Suppose that $c\colon V\tensor V \to V\tensor V$ is a Yang-Baxter operator.
A map $\tilde c\colon V\tensor V \to V\tensor V$ 
is called a \emph{Yang-Baxter deformation} of $c$ (with respect to $\m$)
if $\tilde c$ is itself a Yang-Baxter operator and 
satisfies $\tilde c\equiv c$ modulo $\m$.
\end{definition}

The typical setting is the power series ring $\A = \K\fps{h}$ 
over a field $\K$, equipped with its maximal ideal $\m = (h)$.
In Example \ref{exm:Jones} we can choose $q \in 1+\m$, which ensures 
that $c_q$ is a deformation of $\tau$ in the sense of the definition.

\begin{definition}
An \emph{equivalence transformation} (with respect to the ideal $\m$) 
is an automorphism $\alpha\colon V\to V$ with $\alpha \equiv \I$ modulo $\m$.
Two Yang-Baxter operators $c$ and $\tilde c$ are called 
\emph{equivalent} (with respect to $\m$) if there exists 
an equivalence transformation $\alpha\colon V\to V$ such that 
$\tilde c = (\alpha\tensor\alpha) \, c \, (\alpha\tensor\alpha)^{-1}$.

For every invertible element $s \in 1+\m$ multiplication yields 
a deformation $s\cdot c$ of $c$. Such a rescaling, 
even though uninteresting, is in general not equivalent to $c$. 
A deformation $\tilde c$ of $c$ is called \emph{trivial} 
if it is equivalent to $c$ or to a rescaling $s\cdot c$ 
by some constant factor $s \in 1+\m$. We say that $c$ 
is \emph{rigid} if every Yang-Baxter deformation of $c$ is trivial.
\end{definition}

\subsection{Quandles and racks} \label{sub:Racks}

Racks are a way to construct set-theoretic solutions 
of the Yang-Baxter equation. To begin with, consider a group $G$ 
and a subset $Q \subset G$ that is closed under conjugation.
This allows us to define a binary operation $\ast \colon Q \times Q \to Q$ 
by setting $x\ast y = y^{-1}xy$, which enjoys the following properties:
\begin{enumerate} 
\item[(Q1)]
For every $x\in Q$ we have $x\ast x = x$. 
\hfill (idempotency)
\item[(Q2)]
Every right translation $\varrho(y) \colon x \mapsto x\ast y$ is a bijection.
\hfill (right invertibility)
\item[(Q3)]
For all $x,y,z\in Q$ we have $(x\ast y)\ast z = (x\ast z)\ast (y\ast z)$.
\hfill (self-distributivity)
\end{enumerate}

Such structures have gained independent interest since the 1980s
when they have been applied in low-dimensional topology, 
most notably to study knots and braids.
This is why a general definition has proven useful:

\begin{definition}
Let $Q$ be a set with a binary operation $\ast$.
We call $(Q,\ast)$ a \emph{quandle} if it satisfies axioms (Q1--Q3),
and a \emph{rack} if satisfies axioms (Q2--Q3).
\end{definition}

The term ``quandle'' goes back to D.\,Joyce \cite{Joyce:1982}.
The same structure was called ``distributive groupoid'' 
by S.V.\,Matveev \cite{Matveev:1982}, and 
``crystal'' by L.H.\,Kauffman \cite{Kauffman:2001}.
Since quandles are close to groups, their applications 
in knot theory are in close relationship to the knot group.
We should point out, however, that there exist 
many quandles that do not embed into any group.

Axioms (Q2) and (Q3) are equivalent to saying that every right translation 
$\varrho(y) \colon x \mapsto x\ast y$ is an automorphism of $(Q,\ast)$. 
This is why such a structure was called \emph{automorphic set} 
by E.\,Brieskorn \cite{Brieskorn:1988}. The somewhat shorter term \emph{rack} 
was preferred by R.\,Fenn and C.P.\,Rourke \cite{FennRourke:1992}.

\begin{definition}
Let $Q$ be a rack.
The subgroup of $\Aut(Q)$ generated by the family $\{ \varrho(y) \mid y\in Q \}$ 
is called the group of \emph{inner automorphisms}, denoted $\Inn(Q)$. 
Two elements $x,y\in Q$ are called \emph{behaviourally equivalent},
denoted $x\equiv y$, if $\varrho(x) = \varrho(y)$.
\end{definition}

We adopt the convention that automorphisms of a rack $Q$ 
act on the right, written $x^\phi$ or $x \ast \phi$, 
which means that their composition $\phi\psi$ is defined 
by $x^{(\phi\psi)} = (x^\phi)^\psi$ for all $x\in Q$.
For $x,y \in Q $ we use the notation $x^y$ and $x \ast y$ indifferently.

P.J.\,Freyd and D.N.\,Yetter \cite{FreydYetter:1989}
considered the similar notion of \emph{crossed $G$-sets}.
Here the defining data is a set $Q$ equipped 
with a right action of a group $G$ and 
a map $\varrho\colon Q \to G$ such that
$\varrho(x^g) = g^{-1} \varrho(x) g$. 
One easily verifies that this defines 
a rack $(Q,\ast)$ with $x \ast y = x^{\varrho(y)}$.
Conversely, every rack $(Q,\ast)$ defines a crossed $G$-set 
by choosing the group $G= \Inn(Q)$ with its natural action 
on $Q$ and $\varrho \colon Q \to \Inn(Q)$ as above.
Notice, however, that crossed $G$-sets are slightly more general than quandles,
because the group $G$ acting on $Q$ need not be chosen to be $\Inn(Q)$.

Just as quandles generalize knot colourings, 
racks are tailor-made for braid colourings, 
see E.\,Brieskorn \cite{Brieskorn:1988}.
This brings us back to our main theme:

\begin{proposition}
Given a rack $Q$, one can construct 
a Yang-Baxter operator $c_Q$ as follows: 
let $V = \A{Q}$ be the free $\A$-module with basis $Q$ and define
$$
c_Q \colon \A{Q}\tensor\A{Q} \to \A{Q}\tensor\A{Q} \quad\text{by}\quad
x\tensor y \mapsto y \tensor (x\ast y) \quad\text{for all}\quad x,y\in Q.
$$
By construction, $c_Q$ is a Yang-Baxter operator:
Axiom (Q2) ensures that $c_Q$ is an automorphism,
while Axiom (Q3) implies the Yang-Baxter equation.
\qed
\end{proposition}

\subsection{Entropic maps} \label{sub:EntropicMaps}

In examining deformations of the operator $c_Q$ 
we will encounter certain maps $f \colon \A{Q^n} \to \A{Q^n}$
that respect the inner symmetry of the rack $Q$. 
To formulate this precisely, we introduce some notation.

\begin{definition} \label{def:EntropicMap}
Using graphical notation, a map $f\colon \A{Q^n} \to \A{Q^n}$ 
is called \emph{entropic} with respect to $c_Q$ if it satisfies, 
for each $i=0,\dots,n$, the following equation:
$$
\raisebox{-18pt}{\includegraphics{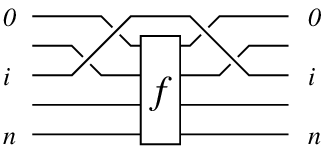}}
\quad = \quad
\raisebox{-18pt}{\includegraphics{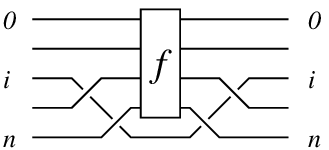}}
$$
\end{definition}
 
This can be reformulated in a more algebraic fashion.
For notational convenience, we do not distinguish between 
the $\A$-linear map $f \colon \A{Q^n} \to \A{Q^n}$ and 
its matrix $f \colon Q^n \times Q^n \to \A$, related by the definition
$$
f\colon \left( x_1 \tensor \cdots \tensor x_n \right) \mapsto
\sum_{y_1,\dots,y_n} f\indices{x_1,\dots,x_n}{y_1,\dots,y_n} 
\cdot \left( y_1 \tensor \cdots \tensor y_n \right) .
$$
Matrix entries are thus denoted by $f\smallindices{x_1,\dots,x_n}{y_1,\dots,y_n}$
with indices $\smallindices{x_1,\dots,x_n}{y_1,\dots,y_n} \in Q^n\times Q^n$.

\begin{definition} \label{def:Quasidiagonal}
A map $f\colon \A{Q^n} \to \A{Q^n}$ is called \emph{quasi-diagonal}
if $f\smallindices{x_1,\dots,x_n}{y_1,\dots,y_n} = 0$
whenever $x_i \not\equiv y_i$ for some index $i \in \{1,\dots,n\}$.
It is \emph{fully equivariant} if it is 
equivariant under the action of $\Inn(Q)^n$, that is
$f\smallindices{x_1,\dots,x_n}{y_1,\dots,y_n} =
\smallindices{x_1\ast{\alpha_1},\dots,x_n\ast{\alpha_n}}%
{y_1\ast{\alpha_1},\dots,y_n\ast{\alpha_n}}$
for all $\alpha_1,\dots,\alpha_n \in \Inn(Q)$.
\end{definition}

\begin{proposition}[proved in \textsection\ref{sub:EntropicCochains}] 
\label{prop:EntropicCharacterization}
An $\A$-linear map $f\colon \A{Q^n} \to \A{Q^n}$ is entropic 
if and only if it is both quasi-diagonal and fully equivariant.
\end{proposition}

\begin{remark}
Entropic maps form a sub-algebra of $\End(\A{Q^n})$.
If $Q$ is trivial, then every map is entropic.
If $\Inn(Q)$ acts transitively on $Q$ and
$\varrho \colon Q \to \Inn(Q)$ is injective,
then the only entropic maps are $\lambda\id$ with $\lambda \in \A$.
There are many examples in between the two extremes. Generally speaking, 
the larger $\Inn(Q)$ is, the fewer entropic maps there are.
\end{remark}

\subsection{Yang-Baxter deformations of racks} \label{sub:Deformations}

As we have seen, each rack $Q$ provides 
a particular solution $c_Q$ to the Yang-Baxter equation.
It appears natural to ask for deformations.
Our main result solves this problem:
under generic hypotheses, every Yang-Baxter deformation
of $c_Q$ is equivalent to an entropic deformation.

\begin{definition}
Every deformation $c$ of $c_Q$ can be written 
as $c = c_Q f$ with $f \equiv \II$ modulo $\m$.
We call such a deformation \emph{entropic} if $f$ is entropic.
\end{definition}

The preliminaries being in place, we can now state the main results:

\begin{theorem}[proved in Section \ref{sec:RackCohomology}]
\label{Thm:InfinitesimalClassification}
Consider the infinitesimal case where $\m^2 = 0$.
Then every entropic deformation of $c_Q$
satisfies the Yang-Baxter equation.
If moreover $|\Inn(Q)|$ is invertible in $\A$,
then every Yang-Baxter deformation of $c_Q$
is equivalent to exactly one entropic deformation.
\end{theorem}

Our approach to prove this theorem is classical:
in the infinitesimal case everything becomes linear in first order terms,
and the Yang-Baxter equation can be recast as a cochain complex.
This can reasonably be called the \emph{Yang-Baxter cohomology}.
It is introduced in Section \ref{sec:YangBaxterCohomology} and
calculated in Section \ref{sec:RackCohomology}.
Having this initial result at hand, we can proceed 
from infinitesimal to complete deformations:

\begin{theorem}[proved in Section \ref{sec:CompleteDeformations}]
\label{Thm:CompleteClassification}
Let $\A$ be a ring that is complete with respect to the ideal $\m$,
and let $Q$ be a rack such that $|\Inn(Q)|$ is invertible in $\A$.
Then every Yang-Baxter deformation of $c_Q \colon \A{Q^2} \to \A{Q^2}$ 
is equivalent to an entropic deformation.
\end{theorem}

Notice that the hypotheses are always satisfied 
for a finite rack $Q$ over the complete local ring 
$\A = \Q\fps{h}$ with its maximal ideal $\m = (h)$.

The preceding theorem ensures that 
we can restrict attention to entropic deformations;
however, not every entropic deformation satisfies the Yang-Baxter equation.
Being entropic suffices in the infinitesimal case, but 
in general higher-order terms introduce further obstructions.
Quite surprisingly, they do not depend on $Q$ at all;
higher-order obstructions are exactly the same as in the quantum case:

\begin{theorem}[proved in Section \ref{sec:HigherObstructions}]
\label{Thm:HigherObstrctions}
Consider a rack $Q$ and its associated Yang-Baxter operator 
$c_Q \colon \A{Q^2} \to \A{Q^2}$ over some ring $\A$.
An entropic deformation $c_Q f$ satisfies the Yang-Baxter equation 
if and only if $\tau f$ satisfies the Yang-Baxter equation,
that is, if and only if $f$ is an r-matrix.
\end{theorem}

The transposition operator $\tau$ does not impose any infinitesimal restrictions; 
the only obstructions are those of degree $2$ and higher.
The preceding theorem says that entropic deformations of $c_Q$
are subject to exactly the same higher-order obstructions 
as deformations of $\tau$, plus the entropy condition 
enforced by a non-trivial inner automorphism group $\Inn(Q)$.
In this sense, entropic Yang-Baxter deformations of $c_Q$ 
are just entropic r-matrices. We have thus reduced the theory of 
formal Yang-Baxter deformations of racks to the quantum case
\cite{Drinfeld:1987,Turaev:1988,Kassel:1995,KasselRossoTuraev:1997}.

\subsection{Applications and examples} \label{sub:Examples}

To simplify notation, we will consider here only
quandles $Q$ that embed into some finite group $G$.
This leads to certain classes of examples where deformations
over $\A = \Q\fps{h}$ are particularly easy to understand.

\begin{remark} 
Consider first a trivial quandle $Q$, with $x\ast y = x$ for all $x,y$,
where $c_Q = \tau$ is simply the transposition operator.
Here our results cannot add anything new:
every map $f \colon \A{Q^n} \to \A{Q^n}$ is entropic,
and so Theorem \ref{Thm:InfinitesimalClassification} simply restates 
that there are no infinitesimal obstructions (every deformation 
of $\tau$ satisfies the Yang-Baxter equation modulo $\m^2$).
There are, however, higher-order obstructions, which we have 
carefully excluded from our discussion: these form a subject 
of their own and belong to the much deeper theory 
of quantum invariants (see Example \ref{exm:Jones}).
\end{remark}

After the trivial quandle, which admits many deformations
but escapes our techniques, let us consider the opposite case 
of a rigid operator:

\begin{corollary}
Let $G$ be a finite centreless group that is 
generated by a conjugacy class $Q$.
Then every Yang-Baxter deformation of $c_Q$ over $\Q\fps{h}$ 
is equivalent to $s\cdot c_Q$ with some constant factor $s \in 1+(h)$.
In other words, $c_Q$ is rigid. 
\qed
\end{corollary}

\begin{example} \label{exm:Dihedral3}
The smallest non-trivial example of a rigid operator 
is given by the set $Q = \{ (12), (13), (23) \}$ 
of transpositions in the symmetric group $S_3$, or equivalently 
the set of reflections in the dihedral group $D_3$.
Ordering the basis $Q\times Q$ lexicographically, 
we can represent $c_Q$ by the matrix
$$
c_Q = \left(\begin{smallmatrix}
 1 & \0 & \0 & \0 & \0 & \0 & \0 & \0 & \0 \\
\0 & \0 & \0 & \0 & \0 & \0 &  1 & \0 & \0 \\
\0 & \0 & \0 &  1 & \0 & \0 & \0 & \0 & \0 \\
\0 & \0 & \0 & \0 & \0 & \0 & \0 &  1 & \0 \\
\0 & \0 & \0 & \0 &  1 & \0 & \0 & \0 & \0 \\
\0 &  1 & \0 & \0 & \0 & \0 & \0 & \0 & \0 \\
\0 & \0 & \0 & \0 & \0 &  1 & \0 & \0 & \0 \\
\0 & \0 &  1 & \0 & \0 & \0 & \0 & \0 & \0 \\
\0 & \0 & \0 & \0 & \0 & \0 & \0 & \0 &  1
\end{smallmatrix}\right).
$$

In the case of the Jones polynomial,
the initial operator $\tau$ is trivial but
the deformation $c_q$ is highly non-trivial.
In the present example, the interesting part 
is the initial operator $c_Q$ itself:
the associated link invariant is the number 
of $3$-colourings, as defined by R.H.\,Fox. 
Unlike $\tau$, the Yang-Baxter operator $c_Q$ does not 
admit any non-trivial deformation over $\Q\fps{h}$.
In this sense it is an isolated solution of the Yang-Baxter equation.
\end{example}

There are also racks in between the two extremes,
which are neither trivial nor rigid. We indicate a class of examples
where every infinitesimal deformation can be integrated,
because higher-order obstructions miraculously vanish.

\begin{corollary} \label{cor:CentralDeformation}
Let $G$ be a finite group, generated by $Q = \cup_i Q_i$, 
where $Q_1,\dots,Q_n$ are distinct conjugacy classes of $G$.
Assume further that the centre $Z$ of $G$ 
satisfies $Z \cdot Q_i = Q_i$ for each $i=1,\dots,n$.
Then every Yang-Baxter deformation of $c_Q$ 
over $\Q\fps{h}$ is equivalent to one of the form
$c(x\tensor y) = s_{ij} \cdot y \tensor x^y$ for $x\in Q_i$ and $y\in Q_j$,
with constant factors $s_{ij} \in 1 + h\Q\fps{h}[Z\times Z]$.
Conversely, every deformation of this form satisfies the Yang-Baxter equation.
\qed
\end{corollary}

\begin{example} \label{exm:Dihedral4}
Consider the set of reflections in the dihedral group $D_4$, that is
$$
Q=\{\; (13)\,,\; (24)\,,\; (12)(34)\,,\; (14)(23) \;\}.
$$
This set is closed under conjugation, hence a quandle.
With respect to the lexicographical basis, $c_Q$ is 
represented by the following permutation matrix:

$$
c_Q = 
\left( \begin{smallmatrix}
 1 & \0 & \0 & \0 & \0 & \0 & \0 & \0 & \0 & \0 & \0 & \0 & \0 & \0 & \0 & \0 \\
\0 & \0 & \0 & \0 &  1 & \0 & \0 & \0 & \0 & \0 & \0 & \0 & \0 & \0 & \0 & \0 \\
\0 & \0 & \0 & \0 & \0 & \0 & \0 & \0 & \0 &  1 & \0 & \0 & \0 & \0 & \0 & \0 \\
\0 & \0 & \0 & \0 & \0 & \0 & \0 & \0 & \0 & \0 & \0 & \0 & \0 &  1 & \0 & \0 \\
\0 &  1 & \0 & \0 & \0 & \0 & \0 & \0 & \0 & \0 & \0 & \0 & \0 & \0 & \0 & \0 \\
\0 & \0 & \0 & \0 & \0 &  1 & \0 & \0 & \0 & \0 & \0 & \0 & \0 & \0 & \0 & \0 \\
\0 & \0 & \0 & \0 & \0 & \0 & \0 & \0 &  1 & \0 & \0 & \0 & \0 & \0 & \0 & \0 \\
\0 & \0 & \0 & \0 & \0 & \0 & \0 & \0 & \0 & \0 & \0 & \0 &  1 & \0 & \0 & \0 \\
\0 & \0 & \0 &  1 & \0 & \0 & \0 & \0 & \0 & \0 & \0 & \0 & \0 & \0 & \0 & \0 \\
\0 & \0 & \0 & \0 & \0 & \0 & \0 &  1 & \0 & \0 & \0 & \0 & \0 & \0 & \0 & \0 \\
\0 & \0 & \0 & \0 & \0 & \0 & \0 & \0 & \0 & \0 &  1 & \0 & \0 & \0 & \0 & \0 \\
\0 & \0 & \0 & \0 & \0 & \0 & \0 & \0 & \0 & \0 & \0 & \0 & \0 & \0 &  1 & \0 \\
\0 & \0 &  1 & \0 & \0 & \0 & \0 & \0 & \0 & \0 & \0 & \0 & \0 & \0 & \0 & \0 \\
\0 & \0 & \0 & \0 & \0 & \0 &  1 & \0 & \0 & \0 & \0 & \0 & \0 & \0 & \0 & \0 \\
\0 & \0 & \0 & \0 & \0 & \0 & \0 & \0 & \0 & \0 & \0 &  1 & \0 & \0 & \0 & \0 \\
\0 & \0 & \0 & \0 & \0 & \0 & \0 & \0 & \0 & \0 & \0 & \0 & \0 & \0 & \0 &  1 \\
\end{smallmatrix} \right) 
$$

By construction, this matrix is a solution of the Yang-Baxter equation.
According to Corollary \ref{cor:CentralDeformation}, 
it admits a $16$-fold deformation $c(\lambda)$ given by

$$
c(\lambda) = c_Q + 
\left( \begin{smallmatrix}
\lambda_1   &\lambda_2   &\0          &\0          &\lambda_3   &\lambda_4   &\0          &\0          &
\0          &\0          &\0          &\0          &\0          &\0          &\0          &\0          \\
\lambda_3   &\lambda_4   &\0          &\0          &\lambda_1   &\lambda_2   &\0          &\0          &
\0          &\0          &\0          &\0          &\0          &\0          &\0          &\0          \\
\0          &\0          &\0          &\0          &\0          &\0          &\0          &\0          &
\lambda_5   &\lambda_6   &\0          &\0          &\lambda_7   &\lambda_8   &\0          &\0          \\
\0          &\0          &\0          &\0          &\0          &\0          &\0          &\0          &
\lambda_7   &\lambda_8   &\0          &\0          &\lambda_5   &\lambda_6   &\0          &\0          \\
\lambda_2   &\lambda_1   &\0          &\0          &\lambda_4   &\lambda_3   &\0          &\0          &
\0          &\0          &\0          &\0          &\0          &\0          &\0          &\0          \\
\lambda_4   &\lambda_3   &\0          &\0          &\lambda_2   &\lambda_1   &\0          &\0          &
\0          &\0          &\0          &\0          &\0          &\0          &\0          &\0          \\
\0          &\0          &\0          &\0          &\0          &\0          &\0          &\0          &
\lambda_6   &\lambda_5   &\0          &\0          &\lambda_8   &\lambda_7   &\0          &\0          \\
\0          &\0          &\0          &\0          &\0          &\0          &\0          &\0          &
\lambda_8   &\lambda_7   &\0          &\0          &\lambda_6   &\lambda_5   &\0          &\0          \\
\0          &\0          &\lambda_9   &\lambda_{10}&\0          &\0          &\lambda_{11}&\lambda_{12}&
\0          &\0          &\0          &\0          &\0          &\0          &\0          &\0          \\
\0          &\0          &\lambda_{11}&\lambda_{12}&\0          &\0          &\lambda_9   &\lambda_{10}&
\0          &\0          &\0          &\0          &\0          &\0          &\0          &\0          \\
\0          &\0          &\0          &\0          &\0          &\0          &\0          &\0          &
\0          &\0          &\lambda_{13}&\lambda_{14}&\0          &\0          &\lambda_{15}&\lambda_{16}\\
\0          &\0          &\0          &\0          &\0          &\0          &\0          &\0          &
\0          &\0          &\lambda_{15}&\lambda_{16}&\0          &\0          &\lambda_{13}&\lambda_{14}\\
\0          &\0          &\lambda_{10}&\lambda_9   &\0          &\0          &\lambda_{12}&\lambda_{11}&
\0          &\0          &\0          &\0          &\0          &\0          &\0          &\0          \\
\0          &\0          &\lambda_{12}&\lambda_{11}&\0          &\0          &\lambda_{10}&\lambda_9   &
\0          &\0          &\0          &\0          &\0          &\0          &\0          &\0          \\
\0          &\0          &\0          &\0          &\0          &\0          &\0          &\0          &
\0          &\0          &\lambda_{14}&\lambda_{13}&\0          &\0          &\lambda_{16}&\lambda_{15}\\
\0          &\0          &\0          &\0          &\0          &\0          &\0          &\0          &
\0          &\0          &\lambda_{16}&\lambda_{15}&\0          &\0          &\lambda_{14}&\lambda_{13}
\end{smallmatrix} \right). 
$$

\smallskip

For every choice of parameters $\lambda_1,\dots,\lambda_{16}$,
the matrix $c(\lambda)$ satisfies the Yang-Baxter equation,
and as a special case we get $c(0)=c_Q$. We finally remark that 
the trace of its square is given by
\begin{align*}
\tr\left[ c(\lambda)^2 \right] 
& = 4(\lambda_{1}+1)^2 + 4\lambda_{4}^2 + 4(\lambda_{13}+1)^2 + 4\lambda_{16}^2 \\
& + 8(\lambda_{6}+1)\lambda_{11} + 8(\lambda_{10}+1)\lambda_{7}
  + 8\lambda_{2}\lambda_{3} + 8\lambda_{14}\lambda_{15}
  + 8\lambda_{5}\lambda_{9} + 8\lambda_{8}\lambda_{12},
\end{align*}
which shows that none of the parameters 
can be eliminated by an equivalence transformation.
This proves anew that the deformed operator $c(\lambda)$ 
is not equivalent to the initial operator $c_Q$.
\end{example}

\begin{remark}
It is amusing to note that the minimal Examples \ref{exm:Jones}, \ref{exm:Dihedral3},
and \ref{exm:Dihedral4} are the first three members of the family formed
by reflections in dihedral groups. The following figure nicely summarizes the point:
\end{remark}

\begin{figure}[hbtp]
  \centering
  \includegraphics{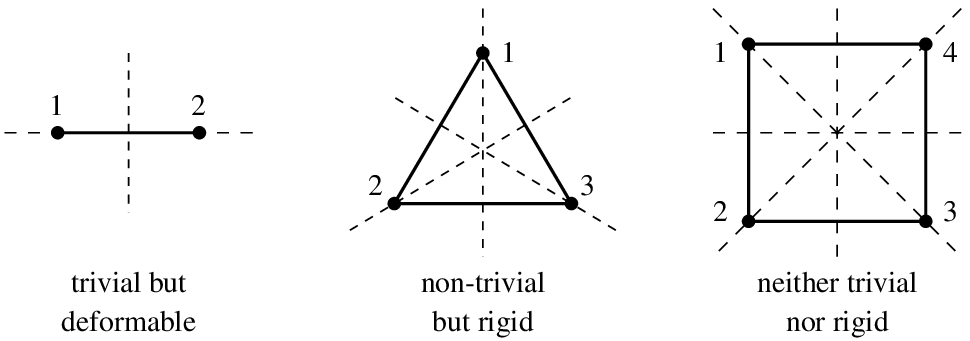}
  \caption{The first three members of the dihedral family}
  \label{fig:DihedralFamily}
\end{figure}


\subsection{Related work} \label{sub:RelatedWork}

Similar deformation and cohomology theories naturally arise 
in situations that are close or equivalent to the Yang-Baxter setting.

\begin{itemize}

\item
Our results can be reformulated in terms of deformations 
of modules over the quantum double $D(G)$ of a finite group $G$.
In this form it has possibly been known to experts in quantum groups,
but there seems to be no written account in the literature.
See \cite[ch.\,IX]{KasselRossoTuraev:1997} for general background. 

\item
The bialgebra approach was pursued by M.\,Gerstenhaber and S.D.\,Schack, 
who  proved in \cite[Section 8]{GerstenhaberSchack:1992} that 
the group bialgebra $\K G$ is rigid as a bialgebra.
They did not discuss deformations of its quantum double $D(G)$.

\item
Our approach can also be reformulated in terms of 
deformations of braided monoidal categories. This point of view 
was put forward by P.J.\,Freyd and D.N.\,Yetter in \cite{FreydYetter:1989}.
The deformation of quandles and racks appeared as an example,
but only diagonal deformations were taken into account.

\item
Diagonal deformations have been more fully developed 
in \cite{CarterEtAl:2003}, where quandle cohomology 
was used to construct state-sum invariants of knots.
P.\,Etingof and M.\,Gra\~na \cite{EtingofGrana:2003} 
have calculated rack cohomology $H^*(Q,\A)$ 
assuming $|\Inn(Q)|$ invertible in $\A$.
Our calculation of $\HYB^2(c_Q,\A)$ 
generalizes their result 
from diagonal to general Yang-Baxter deformations.

\item
In \cite{Yetter:1992} Yetter considered deformations 
of braided monoidal categories in full generality;
see also \cite{Yetter:2001} and the bibliographical references therein.
He was thus led to define a cohomology theory, which is 
essentially equivalent to Yang-Baxter cohomology.
He did not, however, calculate any examples.

\end{itemize}

As far as I can tell, none of the previous results covers 
Yang-Baxter deformations of conjugacy classes, quandles, or racks.


\section{Yang-Baxter cohomology and infinitesimal deformations} \label{sec:YangBaxterCohomology}

This section develops the infinitesimal deformation theory 
of Yang-Baxter operators. As usual, this is most conveniently formulated 
in terms of a suitable cohomology theory, which we will now define.

\subsection{Yang-Baxter cohomology} \label{sub:YangBaxterCohomology}

Let $\A$ be a commutative ring with unit and 
let $\m$ be an ideal in $\A$ . Given an $\A$-module $V$ and 
a Yang-Baxter operator $c \colon V^{\tensor2} \to V^{\tensor2}$, 
we can construct a cochain complex of $\A$-modules
$C^n = \Hom_\A( V^{\tensor n}, \m V^{\tensor n} )$ as follows. 
Firstly, given $f\in C^n$, we define $d^n_i f \in C^{n+1}$ by
$$
d^n_i f = 
\left(c_n \cdots c_{i+1}\right)^{-1} (f \tensor \I) \left(c_n \cdots c_{i+1}\right)
- \left(c_1 \cdots c_i\right)^{-1} (\I \tensor f) \left(c_1 \cdots c_i\right)
$$
or in graphical notation: 
$$
d^n_i f =  
\quad + \quad \raisebox{-18pt}{\includegraphics{coboundary2}}
\quad - \quad \raisebox{-18pt}{\includegraphics{coboundary1}}
$$
We then define the coboundary operator 
$d^n \co C^n \to C^{n+1}$ by $d^n{=}\sum_{i=0}^{i=n} (-1)^i d^n_i$.

\begin{proposition}
The sequence $C^1 \xrightarrow{d^1} C^2 \xrightarrow{d^2} C^3 \dots$
is a cochain complex.
\end{proposition}

\begin{proof}
The hypothesis that $c$ be a Yang-Baxter operator implies 
$d^{n+1}_i d^n_j = d^{n+1}_{j+1} d^n_i$ for $i\le j$.
This can be proven by a straightforward computation;
it is most easily verified using the graphical calculus suggested 
in the above figure. It follows, as usual, that terms cancel 
each other in pairs to yield $d^{n+1} d^n = 0$. 
\end{proof}

\begin{definition}
We call $(C^n,d^n)$ 
the \emph{Yang-Baxter cochain complex} associated with the operator $c$.
As usual, elements of the kernel $Z^n = \Ker(d^n)$ are called \emph{cocycles},
and elements of the image $B^n = \Im(d^{n-1})$ are called \emph{coboundaries}.
The quotient $H^n = Z^n / B^n$ is called the \emph{Yang-Baxter cohomology} 
of the operator $c$, denoted $\HYB^n(c)$, or $\HYB^n(c;\A,\m)$ 
to indicate the dependence on the ring $\A$ and the ideal $\m$.
\end{definition}

\begin{remark}
A more general cohomology can be defined
by taking coefficients in an arbitrary $\A$-module $U$.
The operators $c_i$ act not only on $V^{\tensor n}$ but also on 
$U \tensor V^{\tensor n}$, extended by the trivial action on $U$.
Using this convention, we can define a cochain complex
$C^n = \Hom_\A( V^{\tensor n}, U \tensor V^{\tensor n} )$ 
with coboundary 
given by the same formulae as above. 

Moreover, given a submodule $U' \subset U$, we can consider the image
of the induced map $U' \tensor V^{\tensor n} \to U \tensor V^{\tensor n}$.
(The image will be isomorphic with $U' \tensor V^{\tensor n}$ if $V$ is flat.)
Using this submodule instead of $U \tensor V^{\tensor n}$, 
we obtain yet another cohomology, denoted $\HYB^n(c;U,U')$.
This generalizes our initial definition of $\HYB^n(c;\A,\m)$.
All cohomology calculations in this article 
generalize verbatim to the case $(U,U')$.
For our applications, however, it will be sufficient 
to consider the special case $(\A,\m)$.
\end{remark}

\subsection{Infinitesimal Yang-Baxter deformations} \label{sub:InfinitesimalDeformations}

Consider a Yang-Baxter operator $c \colon V^{\tensor2} \to V^{\tensor2}$.
Every deformation $\tilde c \colon V^{\tensor2} \to V^{\tensor2}$ of $c$ can be written 
as $\tilde c = c(\II+f)$ with perturbation term $f\colon V^{\tensor2} \to \m V^{\tensor2}$.
For the rest of this section we will assume that $\m^2=0$, 
which means that we consider \emph{infinitesimal} deformations.
One can always force this condition by passing to the quotient $\A/\m^2$.
The reason for this simplification is, of course,
that higher-order terms are suppressed and 
everything becomes linear in first order terms.

\begin{proposition}
Suppose that the ideal $\m\subset\A$ satisfies $\m^2 = 0$.
Then $\tilde c = c(\II+f)$ is a Yang-Baxter operator if and only if $d^2 f = 0$.
Moreover, $c$ and $\tilde c$ are equivalent via conjugation by
$\alpha = \I + g$ with $g\colon V \to \m V$ if and only if $f = d^1 g$.
\end{proposition}

\begin{proof}
Spelling out the Yang-Baxter equation for $\tilde c$ yields
the Yang-Baxter equation for $c$ and six error terms of first order.
More precisely, we obtain
\begin{multline*}
(\I\tensor \tilde c)^{-1}(\tilde c \tensor\I)^{-1}(\I\tensor \tilde c)^{-1} 
(\tilde c \tensor\I)(\I\tensor \tilde c)(\tilde c \tensor\I) \\
= 
(\I\tensor c)^{-1}(c \tensor\I)^{-1}(\I\tensor c)^{-1} 
(c \tensor\I)(\I\tensor c)(c \tensor\I) + d^2{f} .
\end{multline*}
By hypothesis, $c$ is a Yang-Baxter operator, 
so the first term is the identity.
As a consequence $\tilde c$ is a Yang-Baxter operator 
if and only if $f \in Z^2(c) := \Ker(d^2)$.

On the other hand, given $\alpha = \I + g$ we have $\alpha^{-1} = \I - g$ and thus
$$
(\alpha\tensor\alpha)^{-1} \, c \, (\alpha\tensor\alpha) = c(\II + d^1{g} )
$$
As a consequence, $c$ and $\tilde c$ are equivalent
if and only if $f \in B^2(c) := \Im(d^1)$.
\end{proof}

The infinitesimal deformations of $c$ 
are thus encoded in the cochain complex
$$
\Hom(V,\m{V}) \xrightarrow{d^1}
\Hom(V^{\tensor 2},\m{V^{\tensor 2}}) \xrightarrow{d^2} 
\Hom(V^{\tensor 3},\m{V^{\tensor 3}}).
$$
Here $d^1$ maps each infinitesimal transformation $g \colon V\to\m{V}$
to its infinitesimal perturbation term $d^1{g} \colon V^{\tensor 2}\to \m{V^{\tensor 2}}$, 
which corresponds to an infinitesimally trivial deformation,
and $d^2$ maps each infinitesimal perturbation $f\colon V^{\tensor 2}\to \m{V^{\tensor 2}}$
to its infinitesimal error term $d^2{f} \colon V^{\tensor 3}\to \m{V^{\tensor 3}}$.
By construction, we find again that $d^2 \compose d^1 = 0$.
We are interested in the quotient $\Ker(d^2) / \Im(d^1)$.



\section{Yang-Baxter cohomology of racks}
\label{sec:RackCohomology}

This section will establish our main technical result:
the explicit calculation of the second Yang-Baxter 
cohomology of a rack $(Q,\ast)$. 
As before, we consider the Yang-Baxter operator 
$c_Q \colon \A{Q^2} \to \A{Q^2}$ defined by 
$x\tensor y \mapsto y \tensor (x\ast y)$.
We wish to study the associated cochain 
complex $C^1 \to C^2 \to C^3 \to \dots$ 
with cocycles $Z^n$ and coboundaries $B^n$.
In degree $2$ this is solved by the following theorem:

\begin{theorem} \label{thm:SecondCohomology}
Entropic $n$-cochains form a submodule of $Z^n$, denoted $E^n$.
If the order of $\Inn(Q)$ is invertible in $\A$,
then we have $Z^2 = E^2 \oplus B^2$, in other words, 
every $2$-cocycle is cohomologous to exactly one entropic cocycle.
\end{theorem}

The theorem implies in particular that $H^2 \cong E^2$,
which is a perfectly explicit description of the second 
Yang-Baxter cohomology of a rack $Q$.
The theorem does even a little better:
in each cohomology class $\xi\in H^2$ 
it designates a preferred representative,
namely the unique entropic cocycle in $\xi$.
This will be proved by a sequence of four lemmas, 
which occupy the rest of this section.

\subsection{The coboundary operators} \label{sub:YBcoboundary}

Our goal is to calculate the Yang-Baxter cohomology of racks.
Before doing so we will first make the coboundary operators more explicit
by translating them from graphical to matrix notation.

Let $\delta \colon Q\times Q \to \A$ be the identity matrix,
which in matrix notation is written as
$$
\delta\indices{x}{y} = \begin{cases}
1 & \text{if } x=y \\ 0 & \text{if } x \ne y \end{cases}.
$$
In this notation the operator
$d_i^n f \colon Q^{n+1} \times Q^{n+1} \to \m$  is given by
\begin{align}
\label{eq:Coboundary}
(d_i^n f)\indices{x_0,\dots,x_n}{y_0,\dots,y_n} = 
+ & f\indices{x_0^{\phantom{x_i}},\dots,x_{i-1}^{\phantom{x_i}},x_{i+1},\dots,x_n}%
{y_0^{\phantom{y_i}},\dots,y_{i-1}^{\phantom{y_i}},y_{i+1},\dots,y_n} 
\cdot \delta\indices{x_i^{x_{i+1}\cdots x_n}}{y_i^{y_{i+1}\cdots y_n}} 
\\ \notag
- & f\indices{x_0^{x_i},\dots,x_{i-1}^{x_i},x_{i+1},\dots,x_n}%
{y_0^{y_i},\dots,y_{i-1}^{y_i},y_{i+1},\dots,y_n} \cdot \delta\indices{x_i}{y_i}.
\end{align}
The coboundary $d^n f \colon C^n \to C^{n+1}$ 
is given by $d^n f = \sum_{i=0}^{i=n} (-1)^i d^n_i f$.

\begin{remark}
Our definitions were motivated by infinitesimal deformations
in the space of Yang-Baxter operators. We could instead restrict 
all coboundary operators to diagonal matrices, that is, to matrices 
$f \colon Q^n \times Q^n \to \m$ with $f\smallindices{x_1,\dots,x_n}{y_1,\dots,y_n} = 0$ 
whenever $x_i \ne y_i$ for some $i$. In this case we obtain the cochain complex 
of quandle or rack cohomology (see 
\cite{CarterEtAl:2003,Eisermann:2003,FennEtAl:2005}).
\end{remark}

\subsection{Characterization of entropic maps} \label{sub:EntropicCochains}

Recall from Definition \ref{def:EntropicMap} 
that a map $f\colon \A{Q^n} \to \m{Q^n}$
is entropic if and only if $d_0 f = \dots = d_n f = 0$.
The following lemma gives a useful reformulation:

\begin{lemma} \label{lem:EntropicCharacterization}
Given an $\A$-linear map $f\colon \A{Q^n} \to \m{Q^n}$ and any $k\in\{0,\dots,n\}$,
we have $d_k f = \dots = d_n f = 0$ if and only if 
the following two conditions hold:
\begin{enumerate}
\item[$D_k$:]
$f\indices{x_1,\dots,x_n}{y_1,\dots,y_n} = 0$
whenever $x_i \not\equiv y_i$ for some $i>k$, and
\item[$E_k$:]
$f\indices{x_1,\dots,x_n}{y_1,\dots,y_n} = 
f\indices{x_1^\alpha,\dots,x_i^\alpha,x_{i+1},\dots,x_n}%
{y_1^\alpha,\dots,y_i^\alpha,y_{i+1},\dots,y_n}$
for all $\alpha\in \Inn(Q)$ and $i\ge k$.
\end{enumerate}
In particular, $f$ is entropic if and only if 
it is quasi-diagonal 
and fully equivariant. 
\end{lemma}

\begin{proof}
By equation \eqref{eq:Coboundary}, conditions 
$D_k$ and $E_k$ imply that $d_k f = \dots = d_n f = 0$.
To prove the converse, we proceed by a downward induction on $k = n,\dots,0$.

Assume $d_k f = \dots = d_n f = 0$ and that $D_{k+1}$ 
and $E_{k+1}$ are true. We want to establish $D_k$ and $E_k$. 
First of all, we can suppose that $x_{k+2} \equiv y_{k+2}$,
\dots, $x_n \equiv y_n$; otherwise $D_k$ and $E_k$ 
are trivially satisfied because all terms vanish.

In order to prove $D_k$, consider the case $x_{k+1} \not\equiv y_{k+1}$.
Since $\varrho(x_{k+1}) \ne \varrho(y_{k+1})$, there exists
$w\in Q$ with such that $u = w \ast {\varrho(x_{k+1})^{-1}}$
differs from $v = w \ast {\varrho(y_{k+1})^{-1}}$.
We can thus choose $u \ne v$ with $u^{x_{k+1}} = v^{y_{k+1}}$
to obtain
$$
0 = (d_k f)\indices{x_1,\dots,x_k,u,x_{k+1},\dots,x_n}{y_1,\dots,y_k,v,y_{k+1},\dots,y_n}
= f\indices{x_1,\dots,x_k,x_{k+1},\dots,x_n}{y_1,\dots,y_k,y_{k+1},\dots,y_n}.
$$

In order to prove $E_k$, it suffices to consider $\alpha = \varrho(z)$ 
with $z\in Q$, since these automorphisms generate $\Inn(Q)$.
Here we obtain
\begin{align*}
0 & = (d_k f)\indices{x_1,\dots,x_k,z,x_{k+1},\dots,x_n}{y_1,\dots,y_k,z,y_{k+1},\dots,y_n} \\
  & = f\indices{x_1,\dots,x_k,x_{k+1},\dots,x_n}{y_1,\dots,y_k,y_{k+1},\dots,y_n}
  - f\indices{x_1^z,\dots,x_k^z,x_{k+1},\dots,x_n}{y_1^z,\dots,y_k^z,y_{k+1},\dots,y_n}.
\end{align*}
This establishes the induction step $k+1 \to k$ and completes the proof.
\end{proof}

Notice that in the preceding lemma we can choose the ideal $\m = \A$;
we thus obtain the characterization of entropic maps announced 
in Proposition \ref{prop:EntropicCharacterization}.

\subsection{Entropic coboundaries vanish}

On our way to establish $Z^2 = E^2 \oplus B^2$,
we are now in position to prove the easy part:

\begin{lemma} \label{lem:Disjoint} 
If the order of the inner automorphism group $G = \Inn(Q)$ 
is not a zero-divisor in $\A$, then $E^n \cap B^n = \{0\}$.
\end{lemma}

\begin{proof}
Consider a coboundary $f = dg$ that is entropic. We have to show that $f=0$.
By the previous lemma, we know that $f$ is quasi-diagonal,
hence we can assume that $x_i \equiv y_i$ for all $i$.
The equation $f = dg$ then simplifies to
\begin{align*}
f\indices{x_1,\dots,x_n}{y_1,\dots,y_n} = 
\sum_{i=1}^{i=n} (-1)^{i-1} \delta\indices{x_i}{y_i}
\biggl( 
  & g\indices{x_1,\dots,x_{i-1},x_{i+1},\dots,x_n}%
{y_1,\dots,y_{i-1},y_{i+1},\dots,y_n} \\
- & g\indices{x_1 \ast x_i,\dots,x_{i-1}\ast x_i,x_{i+1},\dots,x_n}%
{y_1 \ast y_i,\dots,y_{i-1} \ast y_i,y_{i+1},\dots,y_n}
\biggr)
\end{align*}
Using the equivariance under the action of $G^n$, we obtain
\begin{multline*}
|G|^n \cdot f\indices{x_1,\dots,x_n}{y_1,\dots,y_n} = \sum_{\alpha\in G^n} 
f\indices{x_1^{\alpha_1},\dots,x_n^{\alpha_n}}{y_1^{\alpha_1},\dots,y_n^{\alpha_n}} \\
= \sum_{i=1}^{i=n} (-1)^{i-1} \delta\indices{x_i}{y_i} \sum_{\alpha\in G^n} \biggl( g\indices%
{x_1^{\alpha_1},\dots,x_{i-1}^{\alpha_{i-1}},x_{i+1}^{\alpha_{i+1}},\dots,x_n^{\alpha_n}}%
{y_1^{\alpha_1},\dots,y_{i-1}^{\alpha_{i-1}},y_{i+1}^{\alpha_{i+1}},\dots,y_n^{\alpha_n}} \\
- g\indices%
{x_1^{\alpha_1} \ast x_i^{\alpha_i},\dots,x_{i-1}^{\alpha_{i-1}}\ast x_i^{\alpha_i},
x_{i+1}^{\alpha_{i+1}},\dots,x_n^{\alpha_n}}%
{y_1^{\alpha_1} \ast y_i^{\alpha_i},\dots,y_{i-1}^{\alpha_{i-1}} \ast y_i^{\alpha_i},
y_{i+1}^{\alpha_{i+1}},\dots,y_n^{\alpha_n}}
\biggr) 
\end{multline*}
Fix some index $i$ in the outer sum.
We can assume $x_i = y_i$, otherwise 
$\delta\smallindices{x_i}{y_i} = 0$.
Consider further some index $j<i$.
The maps $x_j \mapsto x_j^{\alpha_j} \ast x_i^{\alpha_i}$ and
$y_j \mapsto y_j^{\alpha_j} \ast y_i^{\alpha_i}$ correspond to 
the action of $\alpha_j \alpha_i^{-1} \varrho(x_i) \alpha_i$.
As $\alpha_j$ runs through $G$, the product 
$\alpha_j \alpha_i^{-1} \varrho(x_i) \alpha_i$ also runs through $G$.
This means that in the inner sum over $\alpha\in G^n$, 
all terms cancel each other in pairs.
We conclude that $|G|^n f = 0$, whence $f = 0$.
\end{proof}

\subsection{Making cocycles equivariant by symmetrization}

Given an automorphism $\alpha\in\Aut(Q)$ and a cochain $f\in C^n$, 
we define the cochain $\alpha f \in C^n$ by 
$$
(\alpha f)\indices{x_1,\dots,x_n}{y_1,\dots,y_n} 
:= f\indices{x_1^\alpha,\dots,x_n^\alpha}{y_1^\alpha,\dots,y_n^\alpha}.
$$
It is easily seen that $d(\alpha f) = \alpha(d f)$, hence $\alpha$ 
maps cocycles to cocycles, and coboundaries to coboundaries.
The induced action on cohomology is denoted 
by $\alpha^* \colon \HYB^*(c_Q) \to \HYB^*(c_Q)$.

\begin{lemma} \label{lem:Symmetrization}
Every inner automorphism $\alpha\in\Inn(Q)$ acts trivially on $\HYB^*(c_Q)$.
If the order of the inner automorphism group $G = \Inn(Q)$ is invertible in $\A$, 
then every cocycle is cohomologous to a $G$-equivariant cocycle.
\end{lemma}

\begin{proof}
It suffices to consider inner automorphisms of the form $\alpha = \varrho(z)$ 
with $z\in Q$, since these automorphisms generate $\Inn(Q)$.
For every cocycle $f \in Z^n$ we then have
\begin{multline*}
f\indices{x_1,\dots,x_n}{y_1,\dots,y_n} 
- f\indices{x_1^\alpha,\dots,x_n^\alpha}{y_1^\alpha,\dots,y_n^\alpha}
= (d_n^n f)\indices{x_1,\dots,x_n,z}{y_1,\dots,y_n,z} \\
= (-1)^n \sum_{i=0}^{n-1} (-1)^i (d_i^n f)\indices{x_1,\dots,x_n,z}{y_1,\dots,y_n,z}
= (d^{n-1} g)\indices{x_1,\dots,x_n}{y_1,\dots,y_n}
\end{multline*}
where the cochain $g\in C^{n-1}$ is defined by
$$
g\indices{u_1,\dots,u_{n-1}}{v_1,\dots,v_{n-1}} 
:= (-1)^n f\indices{u_1,\dots,u_{n-1},z}{v_1,\dots,v_{n-1},z}.
$$
This shows that $f - \alpha f = dg$, 
whence $\alpha$ acts trivially on $\HYB^*(c_Q)$.

If the order of $G = \Inn(Q)$ is invertible in $\A$, 
then we can associate to each cochain $f$ a $G$-equivariant
cochain $\bar f = \frac{1}{|G|} \sum_{\alpha\in G} \alpha f$.
If $f$ is a cocycle then so is $\bar f$, and 
both are cohomologous by the preceding argument.
\end{proof}

\subsection{Calculation of the second cohomology group} \label{sub:YBCohomology}

Specializing to degree $2$, the following lemma completes 
the proof of Theorem \ref{thm:SecondCohomology}.

\begin{lemma} \label{lem:SecondCohomology}
Every equivariant $2$-cocycle is cohomologous to an entropic one.
\end{lemma}

\begin{proof}
By hypothesis, we have $d^2 f = 0$, and according 
to Lemma \ref{lem:EntropicCharacterization} 
equivariance 
is equivalent to $d_2^2 f = 0$.
We thus have $d_0^2 f = d_1^2 f$, or more explicitly:
\begin{equation}
\label{eq:Degree2}
f\indices{v,w}{y,z} \left( \delta\indices{u^{vw}}{x^{yz}} - \delta\indices{u}{x} \right)
= f\indices{u,w}{x,z} \delta\indices{v^w}{y^z} - f\indices{u^v,w}{x^y,z} \delta\indices{v}{y}
\end{equation}
for all $u,v,w,x,y,z\in Q$. It suffices to make $f$ quasi-diagonal, that is, 
to ensure $f\smallindices{v,w}{y,z} = 0$ for $v\not\equiv y$ or $w\not\equiv z$.
The left-hand side then vanishes identically, that is $d_0^2 f=0$, 
which entails that the right-hand side also vanishes, whence $d_1^2 f=0$.

First suppose that $w\not\equiv z$. Then there exists 
a pair $(v,y) \in Q\times Q$ with $v\ne y$ but $v^w = y^z$.
If $(u,x)\in Q\times Q$ also satisfies $u\ne x$ and $u^w = x^z$, 
then Equation \eqref{eq:Degree2} implies that
$f\smallindices{v,w}{y,z} = f\smallindices{u,w}{x,z}$.
To see this, notice that $u^w = x^z$ is equivalent to $u^{vw} = x^{yz}$,
because $\varrho(v)\varrho(w) = \varrho(w)\varrho(v^w)$
and $\varrho(y)\varrho(z) = \varrho(z)\varrho(y^z)$,
with $v^w = y^z$ by our assumption.
This allows us to define a $1$-cochain
$$
g\indices{w}{z} = \begin{cases}
0 & \text{if $w \equiv z$, or else} \\
f\smallindices{v,w}{y,z} &\text{with $v\ne y$ such that $v^w = y^z$} .
\end{cases}
$$
According to the preceding argument, $g\smallindices{w}{z}$ 
is independent of the choice of $v,y$. 
In particular $g$ is equivariant since $f$ is.
This implies $d_1^1 g = 0$, hence $d g = d_0^1 g$:
\begin{align*}
(dg)\indices{u,w}{x,z} 
& = g\indices{w}{z} \left( \delta\indices{u^w}{x^z} - \delta\indices{u}{x} \right)
\intertext{This vanishes whenever $w \equiv z$. 
Otherwise we choose $v\ne y$ with $v^w = y^z$ to obtain}
(dg)\indices{u,w}{x,z} 
& = f\indices{v,w}{y,z} \left( \delta\indices{u^w}{x^z} - \delta\indices{u}{x} \right)
= f\indices{v,w}{y,z} \left( \delta\indices{u^{vw}}{x^{yz}} - \delta\indices{u}{x} \right)
\\
& = (d_0^2 f)\indices{u,v,w}{x,y,z} = (d_1^2 f)\indices{u,v,w}{x,y,z} = f\indices{u,w}{x,z}.
\end{align*}
By this construction, $\bar f := f - dg$ is an equivariant cocycle 
satisfying $\bar f\smallindices{u,w}{x,z} = 0$ whenever $w \not\equiv z$.
For $\bar f$ our initial Equation \eqref{eq:Degree2} thus simplifies to
\begin{equation*}
\bar f\indices{v,w}{y,z} \left( \delta\indices{u^v}{x^y} - \delta\indices{u}{x} \right)
= \left( \bar f\indices{u,w}{x,z} - \bar f\indices{u^v,w}{x^y,z} \right) \delta\indices{v}{y} .
\end{equation*}
If $w \equiv z$ but $v \not\equiv y$, then choose $u\ne x$ 
with $u^v = x^y$: the equation reduces to $\bar f\smallindices{v,w}{y,z} = 0$.
This shows that $\bar f$ is quasi-diagonal, in the sense
that $\bar f\smallindices{v,w}{y,z} = 0$ whenever $v\not\equiv y$ or $w\not\equiv z$.
The left-hand side of our equation thus vanishes identically.
The vanishing of the right-hand side is equivalent to
$\bar f\smallindices{u,w}{x,z} = \bar f\smallindices{u\ast\alpha,w}{x\ast\alpha,z}$
for all $\alpha\in\Inn(Q)$. This proves that $\bar f$ is an entropic cocycle, as desired.
\end{proof}

\begin{proof}[Proof of Theorem \ref{thm:SecondCohomology}]
The preceding lemmas allow us to conclude that $Z^2 = E^2 \oplus B^2$,
provided that the order of $G = \Inn(Q)$ is invertible in $\A$.
Firstly, we have $E^n \cap B^n = \{0\}$ by Lemma \ref{lem:Disjoint}.
Moreover, every cocycle is cohomologous to a $G$-equivariant cocycle
by Lemma \ref{lem:Symmetrization}. Finally, in degree $2$ at least, 
every $G$-equivariant cocycle is cohomologous to an entropic cocycle,
by Lemma \ref{lem:SecondCohomology}.
\end{proof}

\begin{question}
Is it true that $Z^n = E^n \oplus B^n$ for all $n>2$ as well?

While all preceding arguments apply to $n$-cochains in arbitrary degree $n$,
the present calculation of $\HYB^2$ seems to work only for $n=2$.
It is quite possible that some clever generalization will work 
for all $n$, but I could not figure out how to do this.
This state of affairs, while not entirely satisfactory, seems acceptable 
because we use only the second cohomology in subsequent applications.
\end{question}


\section{Complete Yang-Baxter deformations} 
\label{sec:CompleteDeformations}

In this section we will pass from infinitesimal to complete deformations.
In order to do so, we will assume that the ring $\A$ is complete 
with respect to the ideal $\m$, that is, we assume that 
the natural map $\A \to \varprojlim \A/\m^n$ is an isomorphism.

\begin{theorem} \label{thm:CompleteDeformations}
Suppose that the ring $\A$ is complete with respect to the ideal $\m$.
Let $Q$ be a rack such that $|\Inn(Q)|$ is invertible in $\A$.
Then every Yang-Baxter operator $c\colon \A{Q^2} \to \A{Q^2}$ with 
$c \equiv c_Q$ modulo $\m$ is equivalent to an entropic deformation of $c_Q$.
More explicitly, there exists $\alpha \equiv \I$ modulo $\m$ such that 
$(\alpha \tensor \alpha)^{-1} \, c \, (\alpha \tensor \alpha) = c_Q f$
with some entropic deformation term $f\colon \A{Q^2} \to \A{Q^2}$,
$f \equiv \II \mod{\m}$.
\end{theorem}

The proof will use the usual induction argument for complete rings.
We will first concentrate on the crucial inductive step:
the passage from $\A/\m^n$ to $\A/\m^{n+1}$.

\subsection{The inductive step}

To simplify notation, we first assume that $\m^{n+1} = 0$.
One can always force this condition 
by passing to the quotient $\A/\m^{n+1}$.

\begin{lemma} \label{lem:InductiveStep}
Consider a ring $\A$ with ideal $\m$ such that $\m^{n+1} = 0$.
Let $c \colon \A{Q^2} \to \A{Q^2}$ be a Yang-Baxter operator that 
satisfies $c \equiv c_Q$ modulo $\m$ and is entropic modulo $\m^n$.
Then $c$ is equivalent to an entropic Yang-Baxter operator. More precisely,
there exists $\alpha\colon \A{Q} \to \A{Q}$ with $\alpha \equiv \I$ modulo $\m^n$, 
such that $(\alpha\tensor\alpha)^{-1} c_n (\alpha\tensor\alpha)$ 
is an entropic deformation of $c_Q$.
\end{lemma}

This lemma obviously includes and generalizes the infinitesimal case $n=1$,
established in Theorem \ref{thm:SecondCohomology}, 
upon which the following proof relies.
The only new difficulty is that higher-order terms of 
the Yang-Baxter equation render the problem non-homogeneous. 
What is left, fortunately, is an affine structure:

\begin{remark}
Suppose $c \colon \A{Q^2} \to \A{Q^2}$ is a Yang-Baxter operator
that satisfies $c \equiv c_Q$ modulo $\m$ and is entropic modulo $\m^n$.
Let $X_n$ be the set of Yang-Baxter operators 
$\tilde c$ 
with $\tilde c \equiv c \mod{\m^n}$.
For each $\tilde c \in X_n$ we have $\tilde c = c(\II+f)$,
with $f \colon \A{Q^2} \to \m^n Q^2$, and it is 
easily verified that $f \in Z^2(c_Q;\A,\m^n)$.
Conversely, every $f \in Z^2(c_Q;\A,\m^n)$
yields a Yang-Baxter operator $c(\II+f) \in X_n$.
In other words, $X_n$ is an affine space over $Z^2(c_Q;\A,\m^n)$.

Likewise, addition of a coboundary $f = dg \in B^2(c_Q;\A,\m^n)$
produces an equivalent deformation $\tilde c = c(\II+f)$. More explicitly,
since $\m^{n+1} = 0$, the map $\alpha = \I+g$ has inverse $\alpha^{-1} = \I-g$;
we thus obtain $\tilde c = (\alpha\tensor\alpha)^{-1} c (\alpha\tensor\alpha)$ as claimed.
\end{remark}

Having the affine structure at hand, 
we can now proceed from $\A/\m^n$ to $\A/\m^{n+1}$:

\begin{proof}[Proof of Lemma \ref{lem:InductiveStep}]
As before let $X_n$ be the set of Yang-Baxter operators 
$\tilde c$ with $\tilde c \equiv c \mod{\m^n}$.
Using the $Z^2$-affine structure on $X_n$, it suffices 
to find \emph{at least one} entropic solution $\tilde c \in X_n$.
Every other solution will then be of the form $\tilde c(\II+f)$
with $f \in Z^2(c_Q;\A,\m^n)$, hence equivalent to $\tilde c(\II+f')$ 
with some entropic $f'$, according to Theorem \ref{thm:SecondCohomology}.
Since the composition of entropic maps is again entropic,
this suffices to prove the lemma.

In order to find an entropic solution $\tilde c \in X_n$, 
we can first of all symmetrize $c$:
given $\alpha\in\Inn(Q)$, we have $\alpha c \in X_n$,
because $c$ is equivariant modulo $\m^n$. This implies 
that $\bar c = \frac{1}{|G|} \sum_{\alpha\in G} \alpha c$
lies in $X_n$, too.
We thus obtain an equivariant operator $\bar c$, 
which we can write $\bar c = c_Q(\II + e)$
with deformation term $e \colon \A{Q^2} \to \m{Q^2}$.
In order to make $e$ quasi-diagonal,
we decompose $e = f + g$ such that 
$$
f\indices{u,v}{x,y} := \begin{cases}
e\smallindices{u,v}{x,y}  & \text{if $u \equiv x$, $v \equiv y$}, \\
0                         & \text{otherwise}
\end{cases} 
$$ 
and 
$$ 
g\indices{u,v}{x,y} := \begin{cases}
0                         & \text{if $u \equiv x$, $v \equiv y$}, \\
e\smallindices{u,v}{x,y}  & \text{otherwise}
\end{cases}
$$
By hypothesis, $e$ is quasi-diagonal modulo $\m^n$, whence we have 
$e\equiv f \mod{\m^n}$ and $g \colon \A{Q^2} \to \m^n Q^2$.
We obtain by this construction a map $\tilde c = c_Q(\II + f)$ 
that is entropic and satisfies $\tilde c \equiv c \mod{\m^n}$.

We claim that $\tilde c$ actually lies in $X_n$,
that is, $\tilde c$ satisfies the Yang-Baxter equation.
To see this, recall that $\bar c$ is a Yang-Baxter operator.
For $\tilde c = \bar c(\II-g)$ we thus obtain
$$
(\I\tensor \tilde c)^{-1}(\tilde c \tensor\I)^{-1}(\I\tensor \tilde c)^{-1} 
(\tilde c \tensor\I)(\I\tensor \tilde c)(\tilde c \tensor\I) = -d^2 g
$$
It is easy to check that the left-hand side is a quasi-diagonal map,
whereas the right-hand side is zero on the quasi-diagonal.
We conclude that \emph{both} must vanish. This means
that $\tilde c$ satisfies the Yang-Baxter equation, as claimed.
\end{proof}

\subsection{From infinitesimal to complete deformations}

To conclude the passage from infinitesimal to complete, 
it only remains to put the ingredients together:

\begin{proof}[Proof of Theorem \ref{thm:CompleteDeformations}]
Starting with $c_1 := c$ for $n=1$, suppose that $c_n = c_Q f_n$
has a deformation term $f_n$ that is entropic modulo $\m^n$.
By Lemma \ref{lem:InductiveStep}, there exists $\alpha_n\colon \A{Q} \to \A{Q}$
with $\alpha_n \equiv \I$ modulo $\m^n$, such that 
$c_{n+1} := (\alpha_n \tensor \alpha_n)^{-1} c_n (\alpha_n \tensor \alpha_n)$
is given by $c_{n+1} = c_Q f_{n+1}$ with $f_{n+1}$ entropic modulo $\m^{n+1}$.
(In fact, the lemma ensures that such a map $\alpha_n$ exists 
modulo $\m^{n+1}$; this can be lifted to a map $\A{Q} \to \A{Q}$, 
which is again invertible because $\A$ is complete.) 
Completeness of $\A$ ensures that we can pass to the limit and define 
the infinite product $\alpha = \alpha_1 \alpha_2 \alpha_3 \cdots$. 
By construction, $(\alpha\tensor\alpha)^{-1} \, c \, (\alpha\tensor\alpha)$ 
is entropic and equivalent to $c$, as desired.
\end{proof}

\section{Entropic deformations and r-matrices}
\label{sec:HigherObstructions}

As we have seen in the preceding theorem,
every Yang-Baxter deformation of $c_Q$ over a complete ring $\A$
is equivalent to an entropic Yang-Baxter deformation.
Conversely, however, not every entropic deformation
gives rise to a Yang-Baxter operator:
being entropic suffices in the infinitesimal case, but 
in general higher-order terms introduce further obstructions. 
Quite surprisingly, they do not depend on $Q$ at all:

\begin{theorem}
Consider a rack $Q$ and its Yang-Baxter operator 
$c_Q \colon \A{Q^2} \to \A{Q^2}$ over some ring $\A$.
An entropic deformation $\tilde c = c_Q f$ satisfies the Yang-Baxter equation 
if and only if $\tilde\tau = \tau f$ satisfies the Yang-Baxter equation,
that is, if and only if $f$ is an r-matrix.
\end{theorem}

As we have seen, the transposition operator $\tau$ does not impose any 
infinitesimal restrictions (the associated coboundary operator vanishes).
The only obstructions are those of higher order. 
The preceding theorem says that the integration of an entropic
infinitesimal deformation to a complete deformation of $c_Q$
entails exactly the same higher-order obstructions as in the quantum case
\cite{Drinfeld:1987,Turaev:1988,Kassel:1995,KasselRossoTuraev:1997}.

\begin{proof}
The theorem follows once we have established the equation
\begin{multline*}
(\I\tensor \tilde c)^{-1} (\tilde c \tensor\I)^{-1} (\I\tensor \tilde c)^{-1} 
(\tilde c \tensor\I) (\I\tensor \tilde c) (\tilde c \tensor\I) \\
= 
(\I\tensor \tilde\tau)^{-1} (\tilde\tau \tensor\I)^{-1} (\I\tensor \tilde\tau)^{-1}
(\tilde\tau \tensor\I) (\I\tensor \tilde\tau) (\tilde\tau \tensor\I) .
\end{multline*}
One way of proving this equality is
by straightforward and tedious calculation.
It seems more convenient, however, to employ 
a suitable graphical notation.
Recall from Definition \ref{def:EntropicMap}
that $f$ being entropic means
\begin{xalignat*}{3}
\raisebox{-10pt}{\includegraphics{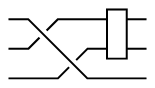}} &= 
\raisebox{-10pt}{\includegraphics{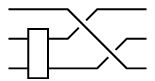}} ,&
\raisebox{-10pt}{\includegraphics{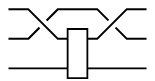}} &= 
\raisebox{-10pt}{\includegraphics{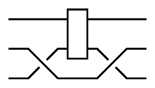}} ,&
\raisebox{-10pt}{\includegraphics{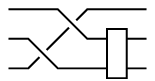}} &= 
\raisebox{-10pt}{\includegraphics{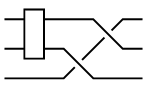}} .
\end{xalignat*}
As before, positive and negative crossings represent $c_Q$ and $c_Q^{-1}$, 
respectively, whereas the box represents the deformation term $f$.
The first and the last equation appear to be rather natural:
they generalize the braid relation (or third Reidemeister move).
The middle equation, however, is somewhat special
and has a curious consequence:

For our operator $c_Q$ we know that the over-passing strand is not affected by a crossing.
The middle equation thus implies that \emph{none} of the strands is affected by 
the shown crossings: we could just as well replace them by transpositions!

Following this observation, our calculation boils down 
to verifying the following transformations:
$$
\begin{matrix}
\raisebox{-10pt}{\includegraphics{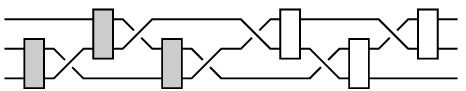}} & = & 
\raisebox{-10pt}{\includegraphics{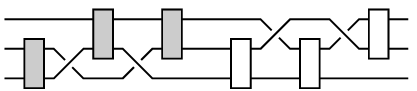}} \\
& & = \\
\raisebox{-10pt}{\includegraphics{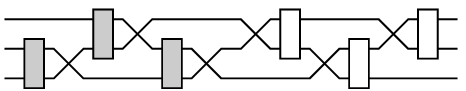}} & = & 
\raisebox{-10pt}{\includegraphics{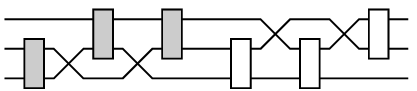}} 
\end{matrix}
$$
Here a white box represents the deformation $f$,
whereas a shaded box represents its inverse $f^{-1}$.
In the first line, positive and negative crossings 
represent $c_Q$ and $c_Q^{-1}$, respectively, whereas 
in the second line, crossings represent the transposition $\tau$.
It is an easy matter to verify the equalities graphically, 
using the fact that $c_Q$ and $\tau$ are Yang-Baxter operators,
and $f$ is entropic with respect to $c_Q$ and $\tau$.
\end{proof}


\section{Closing remarks and open questions} \label{sec:OpenQuestions}


\begin{question}
Our calculation in Section \ref{sec:RackCohomology} relied 
on symmetrization, requiring that $|\Inn(Q)|$ be invertible in $\A$;
this can be seen as the generic case of ``coprime characteristic''.
It seems natural to investigate the cohomology $\HYB^*(c_Q;\A,\m)$ 
in the ``modular case'', where $|\Inn(Q)|$ vanishes in the ring $\A$.
Can one still find a succinct description of $\HYB^2(c_Q)$, 
at least for certain (families of) examples? 
What are the higher order obstructions in this case? 
Do any interesting knot invariants arise in this way?
\end{question}

To illustrate the point, let us emphasize the connection with 
quandle and rack cohomology 
\cite{CarterEtAl:2003,Eisermann:2003,FennEtAl:2005}.
Every map $\alpha \colon Q \times Q \to \Z/n\Z$ defines 
a diagonal deformation of $c_Q$ over $\A = \Z[h]/(n,h^2)$, 
with respect to the ideal $\m = (h)$, by setting 
$c \colon x \tensor y \mapsto \left[ 1 + h\alpha(x,y) \right] \cdot y \tensor x^y$.
One easily checks that $c$ is a Yang-Baxter deformation of $c_Q$ 
if and only if $\alpha$ is an additive rack cocycle, that is
$$
\alpha(x,y) + \alpha(x^y,z) = \alpha(x,z) + \alpha(x^z,y^z). 
$$
Moreover, $c$ is equivalent to $c_Q$ if and only if $\alpha$ is a coboundary.
Sample calculations \cite{CarterEtAl:2003} 
show that rack cohomology $H^2(Q,\Z/n\Z)$ can be non-trivial. P.\,Etingof and M.\,Gra\~na
\cite{EtingofGrana:2003} have shown that this can only happen when $|\Inn(Q)|$ 
is non-invertible in $\A$. In these cases Yang-Baxter cohomology $\HYB^2(c_Q;\A,\m)$ 
will include such extra deformations, and possibly non-diagonal ones, too.



\begin{question}
Is there a topological interpretation of the deformed invariants?
Under suitable conditions, a Yang-Baxter deformation $c$ 
of a quandle $Q$ gives rise to invariants of knots and links
\cite{Turaev:1988,Kassel:1995,KasselRossoTuraev:1997}.
In the case of quandle cohomology one obtains so-called 
state-sum invariants \cite{CarterEtAl:2003}, which have a natural 
interpretation in terms of knot group representations \cite{Eisermann:ColoPoly}.
Can a similar interpretation be established for Yang-Baxter 
deformations of $Q$ in general?
\end{question}

\begin{question}
What can be said about deformations of set-theoretic Yang-Baxter 
operators in general? Following \cite{Drinfeld:1990,ESS:1999,LYZ:2000},
consider a set $Q$ equipped with a bijective map 
$c \colon Q \times Q \to Q \times Q$ satisfying 
$(c \times \I)(\I \times c)(c \times \I) = 
(\I \times c)(c \times \I)(\I \times c)$.
Such a Yang-Baxter map gives rise to a right-action 
$Q \times Q \to Q$, $(x,y) \mapsto x^y$, as well as
a left-action $Q \times Q \to Q$, $(x,y) \mapsto {}^x{\!}y$,
defined by $c(x,y) = ({}^x{\!}y,x^y)$. Notice that the case of 
a trivial left action, ${}^x{\!}y = y$, corresponds precisely to racks.


As for racks, we can extend $c$ to an $\A$-linear map 
$c \colon \A{Q^2} \to \A{Q^2}$ and study deformations over $(\A,\m)$.
As before, Yang-Baxter cohomology $\HYB^2(c_Q;\A,\m)$ yields a convenient framework
for the infinitesimal deformations of $c$, but concrete calculations 
are more involved. Does the cohomology $\HYB^2(c)$ still 
correspond to entropic maps (under suitable hypotheses)?
Can one establish similar rigidity properties?
What happens in the modular case? 
Is there a topological interpretation?
\end{question}


\Addresses\recd


\end{document}
